\documentclass {article}
\usepackage{authblk}
\usepackage{etex}
\usepackage[utf8]{inputenc}
\usepackage[english]{babel}
\usepackage{amsmath}
\usepackage{amsthm}
\usepackage{amssymb}
\usepackage{sectsty}
\usepackage{titlesec}
\usepackage{color}
\usepackage[color,matrix,arrow]{xy}
\usepackage{amsgen}
\usepackage{amstext}
\usepackage{amsbsy}
\usepackage{amsopn}
\usepackage{amsfonts}
\usepackage{eepic}
\usepackage{graphicx}
\usepackage{epsf}
\usepackage{pstricks}
\usepackage{enumerate}
\usepackage{eqnarray}
\usepackage{faktor}
\xyoption{all}

\usepackage{pgf}
\usepackage{tikz}
\usepackage{tikz-cd}

\usetikzlibrary{arrows}

\newcommand{\footrecall}[1]{%
}

\titleformat*{\section}{\large\bfseries}
\titleformat*{\subsection}{\normalsize \bfseries}

\newcommand{\N}{\mathbb{N}}
\newcommand{\Z}{\mathbb{Z}}

\newcommand{\Fix}{\text{Fix}}
\newcommand{\rank}{\text{rank}}

\newcommand{\End}{\text{End}}
\newcommand{\Ker}{\text{Ker}}
\newcommand{\Evper}{\text{EvPer}}
\newcommand{\Evfix}{\text{EvFix}}
\newcommand{\Per}{\text{Per}}

\newcommand{\Orb}{\text{Orb}}
\newcommand{\Aut}{\text{Aut}}
\newcommand{\Img}{\text{Im}}
\newcommand{\mc}{\mathcal}

\theoremstyle{definition}
\newtheorem{theorem}{Theorem}[section]
\newtheorem{corollary}[theorem]{Corollary}

\newtheorem{proposition}[theorem]{Proposition}
\newtheorem{lemma}[theorem]{Lemma}
\newtheorem{example}[theorem]{Example}
\newtheorem{remark}[theorem]{Remark}
\newtheorem{problem}[theorem]{Problem}

\begin{document}
 

\title{Eventually fixed points of endomorphisms of virtually free groups}
\author{Andr\'e Carvalho} 
\maketitle

\begin{abstract}
We consider the subgroup of points of finite orbit through the action of an endomorphism of a finitely generated virtually free group, with particular emphasis on the subgroup of eventually fixed points, $\Evfix(\varphi)$: points whose orbit contains a fixed point. We provide an algorithm to compute the subgroup of fixed points of an endomorphism of a finitely generated virtually free group and prove that finite orbits have cardinality bounded by a computable constant, which allows us to solve several algorithmic problems: deciding if $\varphi$ is a finite order element of $\End(G)$,  if $\varphi$ is aperiodic, if $\Evfix(\varphi)$ is finitely generated and, in the free group case, whether $\Evfix(\varphi)$ is a normal subgroup of $F_n$ or not. We also present a bound for the rank of $\Evfix(\varphi)$ in case it is finitely generated.
\end{abstract}
%
%
\section{Introduction}
The study of fixed subgroups of endomorphisms of groups started with the (independent) work of Gersten \cite{[Ger87]} and Cooper \cite{[Coo87]}, using respectively graph-theoretic
and topological approaches. They proved that the subgroup of fixed points $\Fix(\varphi)$ of
some fixed automorphism $\varphi$ of $F_n$ is always finitely generated, and Cooper succeeded
on classifying from the dynamical viewpoint the fixed points of the continuous extension of $\varphi$ to the boundary of $F_n$. Bestvina and Handel subsequently developed
the theory of train tracks to prove that $\Fix(\varphi)$ has rank at most $n$ in $\cite{[BH92]}$. This was shown to hold for general endomorphisms by Imrich and Turner by reducing to problem to the automorphism case in \cite{[IT89]}. The problem of computing a basis for  $\Fix(\varphi)$  had a tribulated history and was finally settled
by Bogopolski and Maslakova in 2016 in \cite{[BM16]} for automorphisms and by Mutanguha \cite{[Mut21]} for general endomorphisms of a free group. This line of research extended early to wider classes of groups. For instance, Paulin proved in 1989 that the subgroup
of fixed points of an automorphism of a hyperbolic group is finitely generated \cite{[Pau89]}.
Fixed points were also studied for right-angled Artin groups \cite{[RSS13]} and lamplighter
groups \cite{[MS21]}.

In this paper, we extend the algorithmic result of \cite{[Mut21]} proving that $\Fix(\varphi)$ is computable if $\varphi$ is an endomorphism of a finitely generated virtually free group.
\newtheorem*{computvfree}{Theorem \ref{computvfree}}
\begin{computvfree}
Let $G$ be a finitely generated virtually free group and $\varphi\in \End(G)$. Then $\Fix(\varphi)$ is computable. 
\end{computvfree}

In \cite{[MS02]}, Myasnikov and Shpilrain study finite orbits of elements of a free group under the action of an automorphism proving that, in a free group $F_n$, there is an orbit of cardinality $k$ if and only if there is an element of order $k$ in $\Aut(F_n)$. Moreover, the authors prove that this result does not hold for general endomorphisms, by providing an example of an endomorphism of $F_3$ for which there is a point whose orbit has $5$ elements.

In this paper, we study finite orbits of elements under the action of an endomorphism of a finitely generated virtually free group. The orbit of an element $x\in F_n$ through an endomorphism $\varphi$ is finite if and only if it intersects the subgroup of periodic points, $\Per(\varphi)$, of $\varphi$. The set of such points forms a subgroup of $F_n$ and so do the points whose orbit intersects the fixed subgroup $\Fix(\varphi)$. We call these subgroups $\Evper(\varphi)$ and $\Evfix(\varphi)$, respectively. It is easy to see that $\Evfix(\varphi)$ coincides with $\Fix(\varphi)$ if (and only if) $\varphi$ is injective. For this reason, we will mainly focus on noninjective endomorphisms. Given an endomorphism $\varphi\in \End(F_n)$, we will find conditions for $\Evfix(\varphi)$ to be a normal subgroup of $F_n$. Also, despite the fact that the result in \cite{[MS02]} cannot be generalized to endomorphisms, we prove that, replacing finite orbits by periodic orbits, the result holds for endomorphisms, i.e.,  there is an endomorphism of $F_n$ with a periodic orbit of cardinality $k$ if and only if there is an element of order $k$ in $\Aut(F_n)$. Moreover, we prove that, for endomorphisms of finitely generated virtually free groups, finite orbits have bounded cardinality.

\newtheorem*{upbound}{Corollary \ref{upbound}}
\begin{upbound}
Let $G$ be a finitely generated virtually free group and $\varphi\in \End(G)$. There is a computable constant $k$ that 
$$\max\{|\Orb_\varphi(x)| \, \big\lvert \, x\in \Evper(\varphi)\}\leq k.$$
\end{upbound}

This allows us to solve some algorithmic questions: we can decide if $\varphi$ is a finite order element of $\End(G)$, if $\varphi$ is aperiodic or not and, in case $G$ is free, whether $\Evfix(\varphi)$ is a normal subgroup of $F_n$ or not.

Unlike the case of $\Fix(\varphi)$ and $\Per(\varphi)$, the subgroups $\Evfix(\varphi)$ and $\Evper(\varphi)$ are not necessarily finitely generated. However, we prove that we can always decide if that is the case, by proving the following result, which might have independent interest. 
Also, if $\Evfix(\varphi)$ is finitely generated, then a set of generators can be effectively computed.

\newtheorem*{decidepreimage}{Proposition \ref{decidepreimage}}
\begin{decidepreimage}
Let $G$ be a finitely generated virtually free group having a free subgroup $F$ of finite index, $H\leq G$ and $\varphi\in \End(G)$ be an endomorphism. If $\Ker(\varphi)$ is finite, then $H\varphi^{-1}$ is finitely generated. If not, the following are equivalent:
\begin{enumerate}
\item $H\varphi^{-1}$ is finitely generated
\item $H\varphi^{-1}\cap F$ is a finite index subgroup of $F$
\item  $H\varphi^{-1}\cap F$ is a finite index subgroup of $G$
\item $H\varphi^{-1}$ is a finite index subgroup of $G$
\item $H\cap G\varphi$ is a finite index subgroup of $G\varphi$
\item $H\cap F\varphi$ is a finite index subgroup of $F\varphi$
\end{enumerate}
\end{decidepreimage}

Finally, we provide an upper bound for the rank of $\Evfix(\varphi)$ for endomorphisms of $G$. However, we do not know if the bound is sharp.


\section{General properties}
Let $G$ be a group and $\varphi\in \End(G)$. A point $x\in G$ is said to be a \emph{fixed point} if $x\varphi=x$. The set of all fixed points forms a subgroup which we denote by $\Fix(\varphi)$. A point  $x\in G$ is said to be a \emph{periodic point} if there is some $m\in \N$ such that  $x\varphi^m=x$. The set of all periodic points forms a subgroup which we denote by $\Per(\varphi)$. Obviously, we have that $$\Per(\varphi)=\bigcup_{k= 1}^\infty \Fix(\varphi^k).$$

Given $x\in G$, the \emph{orbit of $x$ through $\varphi$} is defined by $$\Orb_\varphi(x)=\{x\varphi^k\mid k\in \N\}.$$ A point $x$ is said to be \emph{eventually periodic} if its orbit is finite, i.e., there is some $m\in \N$ such that $x\varphi^m\in \Per(\varphi)$  and similarly, $x$ is said to be  \emph{eventually fixed} if there is some $m\in \N$ such that $x\varphi^m\in \Fix(\varphi)$ or, equivalently, such that $x\varphi^m=x\varphi^{m+1}$. In this case, for every $k\geq m$, we have that $x\varphi^k=x\varphi^m$. We denote by $\Evper(\varphi)$ (resp. $\Evfix(\varphi)$) the set of all eventually periodic (resp. fixed) points of $\varphi$. It is clear from the definitions that $$\Evper(\varphi)=\bigcup_{k=1}^\infty (\Per(\varphi))\varphi^{-k}\quad \text{and} \quad \Evfix(\varphi)=\bigcup_{k=1}^\infty (\Fix(\varphi))\varphi^{-k}.$$

\begin{proposition}
Let $\varphi\in \End(G).$ Then $\Evper(\varphi)$ and $\Evfix(\varphi)$ are groups.
\end{proposition}
\noindent\textit{Proof.} Let  $x_1,x_2\in \Evper(\varphi)$. Then, there are $m_1,m_2\in \N$ such that  $x_1\varphi^{k_1}$ and $x_1\varphi^{k_2}$ are periodic points for all $k_1>m_1$ and $k_2>m_2$. So, taking $M=\max\{m_1,m_2\}$, we have that $(x_1x_2)\varphi^M$ is periodic. Also, if there is some $m\in \N$ such that $x_1\varphi^m\in \Per(\varphi)$, then $x_1^{-1}\varphi^m\in \Per(\varphi).$

 Similarly, let $x_1,x_2\in \Evfix(\varphi)$. Then, there are $m_i\in \N$ such that  $x_i\varphi^{m_i}=x_i\varphi^{m_i+1}$, for $i=1,2$. Then, putting $M=\max\{m_1,m_2\}$, we have that  $(xy)\varphi^{M+1}=x\varphi^{M+1}y\varphi^{M+1}=x\varphi^{m_1}y\varphi^{m_2}=x\varphi^{M}y\varphi^{M}.$  Also, we have that $(x_1^{-1})\varphi^{m_1+1}=(x_1\varphi^{m_1+1})^{-1}=(x_1\varphi^m)^{-1}=x_1^{-1}\varphi^m$. 
\qed

Now, we present some natural properties of $\Evfix(\varphi)$.
\begin{lemma}
Let $\varphi\in \End(G).$ Then
\begin{enumerate}[(i)]
\item $\bigcup_{k=1}^\infty \Ker(\varphi^k)\trianglelefteq \Evfix(\varphi)$
\item $\Fix(\varphi)\leq \Evfix(\varphi)$
\item $\Evfix(\varphi)\cap \Per(\varphi)=\Fix(\varphi)$
\item $\Evfix(\varphi)=\Fix(\varphi) \Leftrightarrow \varphi \text{ is a monomorphism}$
\end{enumerate}
\end{lemma}
\noindent\textit{Proof.} (i), (ii) and (iii) are obvious by definition. If $\varphi$ is a monomorphism, then $(\Fix(\varphi))\varphi^{-k}=\Fix(\varphi)$ for every $k\in \N$, and so $\Evfix(\varphi)=\Fix(\varphi).$ If there exists $1\neq x\in \Ker(\varphi),$ then $x\in \Evfix(\varphi)\setminus \Fix(\varphi).$
\qed\\

In case $G$ is a finitely generated virtually free group, it is well known that $\Fix(\varphi)$ is finitely generated. However, that might not be the case of $\Evfix(\varphi)$, even in the free group case. 
 
\begin{example}
\label{exnotfg}
Let $\varphi:F_2\to F_2$ be defined by $a\mapsto aba$ and $b\mapsto 1$. Then for $x\in F_2$, we have that $x\varphi=(aba)^{\lambda_a(x)}$, where $\lambda_a:F_2\to \Z$ is the endomorphism defined by $a\mapsto 1$ and $b\mapsto 0.$ So, $\Fix(\varphi)$ is trivial and $\Ker(\varphi)=\{w\,\lvert \, \lambda_a(w)=0\}$. Also $(\Ker(\varphi))\varphi^{-1}=\Ker(\varphi)$. So, $\Evfix(\varphi)=\Ker(\varphi)$, which is not finitely generated.
\end{example}

\begin{example}
\label{exfg}
Let $\varphi:F_2\to F_2$ be defined by $a\mapsto bab^{-1}$ and $b\mapsto 1$. Then for $x\in F_2$, we have that $x\varphi=ba^{\lambda_a(x)}b^{-1}$. So, $\Fix(\varphi)=\{ba^kb^{-1}\,\lvert\, k\in \mathbb Z\}$ and $\Evfix(\varphi)=(\Fix(\varphi))\varphi^{-1}=F_2$, which is finitely generated.
\end{example}

We will see later in the paper that we can decide whether $\Evfix(\varphi)$ and $\Evper(\varphi)$ are finitely generated for endomorphisms of finitely generated virtually free groups.


\section{Normality in free groups}

The purpose of this section is to describe the cases where, for an endomorphism $\varphi\in \End(F_n)$, we have that $\Evfix(\varphi)\trianglelefteq F_n$. We start with a technical lemma.
\begin{lemma}
\label{conjfree}
Let $u \in F_n \setminus \{1\}$
be a nontrivial non proper power. If there exists $w \in F_n$ and $p,q \in \Z$ such that $w^{-1} u^p w = u^q$, then $p = q$ and $w\in \langle u \rangle$.
\end{lemma}

\noindent\textit{Proof.} Let $u \in F_n \setminus \{1\}$. If there exist $w \in F_n$ and $p,q \in \Z$ such that $w^{-1} u^p w = u^q$, then 
the cyclic reduced cores of $u^p$ and $u^q$ are equivalent under a cyclic permutation of their letters. This implies in particular that $p=q$ and $w$ commutes with $u^p.$
\qed\\

We are now able to describe the cases where $\Evfix(\varphi)$ is a normal subgroup of $F_n$.
\begin{proposition}
\label{normalevfix}
Let $\varphi\in\End(F_n).$ Then one of the following holds:
 \begin{enumerate}
 \item $\Evfix(\varphi)=F_n$
 \item $\Evfix(\varphi)=\bigcup_{k=1}^\infty \Ker(\varphi^k)$
 \item $\Evfix(\varphi)$ is not a normal subgroup of $F_n$.
 \end{enumerate}
\end{proposition}
\noindent\textit{Proof.} Suppose that  $\Evfix(\varphi)\neq F_n$, $\Evfix(\varphi)\neq \bigcup_{k=1}^\infty \Ker(\varphi^k)$ 
and that $\Evfix(\varphi)\trianglelefteq F_n$ and take $g\in F_n\setminus \Evfix(\varphi)$. Since $ \Evfix(\varphi)$ is normal, then for every
 $x\in \Evfix(\varphi)\setminus \bigcup_{k=1}^\infty \Ker(\varphi^k)$
we have that  $gxg^{-1}\in \Evfix(\varphi)$, and $g^{-1}xg\in \Evfix(\varphi)$, which means that there are some $n,m\in \N$ for which $(gxg^{-1})\varphi^{n+1}=(gxg^{-1})\varphi^{n}$ and $(g^{-1}xg)\varphi^{m+1}=(g^{-1}xg)\varphi^{m}$. 
Also, there is some $p\in \N$ such that $x\varphi^{p}=x\varphi^{p+1}$. 
So, letting $M=\max\{n,m,p\}$, it follows that 
\begin{align}
\label{gxg-1}
(gxg^{-1})\varphi^{M+1}=(gxg^{-1})\varphi^{M},
\end{align}
 \begin{align}
 \label{g-1xg}
 (g^{-1}xg)\varphi^{M+1}=(g^{-1}xg)\varphi^{M}
 \end{align}
  and 
  \begin{align}
  \label{x}
  x\varphi^{M}=x\varphi^{M+1}.
  \end{align}
 We can rewrite (\ref{gxg-1}) as
$$x\varphi^{M+1}g^{-1}\varphi^{M+1}g\varphi^{M}=g^{-1}\varphi^{M+1}g\varphi^{M}x\varphi^{M},$$ 
and so $g^{-1}\varphi^{M+1}g\varphi^{M}$ and $x\varphi^{M}$ commute. 
Similarly, we can rewrite (\ref{g-1xg}) as
$$g\varphi^{M}g^{-1}\varphi^{M+1}x\varphi^{M+1}=x\varphi^{M}g\varphi^{M}g^{-1}\varphi^{M+1},$$

Since  $g\not\in \Evfix(\varphi)$, then $g\varphi^{M}g^{-1}\varphi^{M+1},g^{-1}\varphi^{M+1}g\varphi^{M}\neq 1$, and since $x\not\in \bigcup_{k=1}^\infty \Ker(\varphi^k)$, then $x\varphi^{M}\neq 1$.
We then have that $g\varphi^{M}g^{-1}\varphi^{M+1},g^{-1}\varphi^{M+1}g\varphi^{M}$ and $x\varphi^{M}$  are powers of the same primitive word $u\in F_n$. So, put 
\begin{align}
\label{c1}
g^{-1}\varphi^{M+1}g\varphi^{M}=u^{q}
\end{align}
\begin{align} 
\label{c2}
g\varphi^{M}g^{-1}\varphi^{M+1}=u^{k}
\end{align} 
and 
\begin{align}
\label{c3}
x\varphi^{M}=u^{r}.
\end{align} 

From  (\ref{c1}) we get that $$g^{-1}\varphi^{M+1}g\varphi^Mg^{-1}\varphi^{M+1}g\varphi^M=u^{2q}.$$ Applying (\ref{c2}), we obtain that $g^{-1}\varphi^{M+1}u^kg\varphi^M=u^{2q},$ and so that
$$g\varphi^{M+1}=u^kg\varphi^Mu^{-2q}.$$ But from  (\ref{c1}), we know that $$g\varphi^{M+1}=g\varphi^{M}u^{-q}.$$

Hence $g\varphi^{M}u^{-q}=u^kg\varphi^Mu^{-2q}$, and $g\varphi^{M}u^{q}g^{-1}\varphi^{M}=u^k.$  From Lemma \ref{conjfree}, we know that $k=q$ and $g\varphi^{M}$ is a power of $u$.

Since $x\varphi^M\neq 1$, then $r\neq 0$. From (\ref{x}) and (\ref{c3}), we know that $(u\varphi)^r=u^r\varphi=u^r$ and so $u\in \Fix(\varphi)$ 
and $g\varphi^{M}$ is fixed, which contradicts the assumption that $g\not\in \Evfix(\varphi).$

\qed\\

We now present some remarks on the previous Proposition.

\begin{remark}
\label{c2}
Condition 2. in Proposition \ref{normalevfix} is equivalent to $\Fix(\varphi)$ being trivial. Indeed, if there is some nontrivial element $x\in \Fix(\varphi)$, then $x\varphi^k=x\neq 1$, for every $k\in \N$. If $ \Fix(\varphi)=1$, then 2. holds by definition. 
\end{remark}

\begin{remark}
If $\Evfix(\varphi)$ is finitely generated, then there must be a bound on the size of the orbits of eventually fixed points, but the converse is false by Example \ref{exnotfg}. Indeed, suppose that $\Evfix(\varphi)=\langle w_1,\ldots, w_k\rangle$. Then 
$$M=\max\{|\Orb_\varphi(w_i)| \, \big \lvert\, i\in [k]\}$$
is a bound on the size of the finite orbits. We will see later that such a bound always exists if $\varphi$ is an endomorphism of a finitely generated virtually free group.
\end{remark}
\begin{remark}
Conditions 1. and 2. are not mutually exclusive. However, if $F_n=\bigcup_{k=1}^\infty \Ker(\varphi^k)$, then every point is eventually sent to $1$ and, by the observation above, orbit sizes must be bounded. So, this happens if and only if $\varphi$ is a \emph{vanishing endomorphism}, i.e, if there is some $r$ such that 
$\varphi^r$ maps every element to $1$. 
\end{remark}


\section{An algorithm to compute the fixed subgroup of an endomorphism of a finitely generated virtually free group}

The purpose of this section is to provide an algorithm to compute the fixed subgroup of an endomorphism of a finitely generated virtually free group. 
When we take a finitely generated virtually free group as input, we assume that we are given a decomposition as a disjoint union
\begin{align}
\label{decomp}
G=Fb_1\cup Fb_2\cup \cdots \cup Fb_m,
\end{align}
where $F=F_A\trianglelefteq G$ is a finitely generated free group and a presentation of the form $\langle A,b_1,\ldots, b_m \mid R\rangle$, where the relations in $R$ are of the  form $b_ia=u_{ia}b_{i}$ and $b_{i}b_{j}=v_{ij}b_{r_{ij}}$, with  $u_{ia}, v_{ij} \in F_A$ and $r_{ij}\in [m]$, $i,j=1\ldots,m$, $a\in A$.

We start by presenting a technical lemma, which is simply an adaptation of \cite[Lemma 2.2]{[HW03]} with the additional condition of the subgroups being normal. The proof follows in the exact same way as theirs, noting that the preimage of a normal subgroup by an endomorphism is still a normal subgroup.

A subgroup $H$ of a group $G$ is \emph{fully invariant} if $\varphi(H)\subseteq H$ for every endomorphism $\varphi$ of $G$.

\begin{lemma}
\label{hsuwise}
Let $G$ be a group, $n$ be a natural number and $N$ be the intersection of all normal subgroups of $G$ of index $\leq n$. Then $N$ is fully invariant, and if $G$ is finitely generated, then $N$ has finite index in $G$.
\end{lemma}

\begin{theorem}
\label{computvfree}
Let $G$ be a finitely generated virtually free group and $\varphi\in \End(G)$. Then $\Fix(\varphi)$ is computable. 
\end{theorem}

\noindent\textit{Proof.}  Take a decomposition as in (\ref{decomp}).
By Lemma \ref{hsuwise}, the intersection $F'$ of all  normal subgroups of $G$ of index at most $m$  is a fully invariant finite index subgroup, i.e. $[G:F']<\infty$ and $F'\varphi\subseteq F'$ for all endomorphisms $\varphi\in \End(G)$. Also, since $[G:F]=m$ and $F\trianglelefteq G$, then $F'\leq F$, and so $F'$ is free. We will now prove that $F'$ is computable. We start by  enumerating all finite groups of cardinality at most $m$. For each such group $K=\{k_1,\ldots, k_s\}$ we enumerate all homomorphisms from $G$ to $K$ by defining images of the generators and checking all the relations. For each homomorphism $\theta:G\to K$, we have that $[G:\Ker(\theta)]=|\Img(\theta)|\leq |K|\leq m$. In fact, all normal subgroups of $G$ of index at most $m$ are of this form. We compute generators for the kernel of each $\theta$, which is possible since we can test membership in $\Ker(\theta)$, which is a finite index subgroup. We can also find $a_2,\ldots,a_{s}\in G$ such that 
$$G=\Ker(\theta)\cup\Ker(\theta)a_2\cup\cdots\cup \Ker(\theta)a_{s},$$
taking $a_i$ such that $a_i\theta=k_i.$

Hence, we can compute $F'$, since it is a finite intersection of computable subgroups, and a decomposition of $G$ as a disjoint union

\begin{align}
\label{decomp2}
G=F'b_1'\cup Fb_2'\cup \cdots \cup Fb_m'.
\end{align}

Now, take $\psi=\varphi|_{F'}$. Since $F'$ is fully invariant, then $\psi\in \End(F')$.  Now, for $u\in F'$, put $X_u=\{x\in F'\mid x\psi=xu\}$. We claim that $X_u$ is computable. Since $[G:F']<\infty$, then $F'$ is finitely generated and so it has a finite basis $X$.
Consider a new letter $c$ not belonging $X$, let $F''=F'*\langle c| \rangle$ and $\psi'\in \End(F'')$ defined by mapping the letters $x\in F'$ to $x\psi$ and $c$ to $u^{-1}c$. By \cite{[Mut21]}, we can compute a basis for $\Fix(\psi')$. It is easy to see that $X_uc=\Fix(\psi')\cap F'c$. Indeed, if $x\in X_u$, then $$(xc)\psi'=(x\psi')(c\psi')=(x\psi)u^{-1}c=xuu^{-1}c=xc$$ and if $x\in \Fix(\psi')\cap F'c$, then there is $y\in F'$ such that $x=yc$ and $$yc=(yc)\psi'=(y\psi)u^{-1}c,$$ which means that $y\psi=yu$ and so $y\in X_u$. Therefore, for all $u\in F'$, $X_uc$ (and so $X_u$) is computable.
We claim that, for $i\in [m']$, 
 $$\Fix(\varphi)\cap F'b_i'=
 \begin{cases}
  \emptyset \quad&\text{if $b_i'((b_i')^{-1}\varphi)\not\in F'$}\\
  X_{b_i'((b_i')^{-1}\varphi)}b_i' \quad&\text{if $b_i'((b_i')^{-1}\varphi)\in F'$}
 \end{cases}.
 $$ 
  Suppose that $b_i'((b_i')^{-1}\varphi)\in F'$. Let $x\in  X_{b_i'((b_i')^{-1}\varphi)}$. Then, $x\psi=xb_i'((b_i')^{-1}\varphi)$, and so $(xb_i')\varphi=xb_i'$. Thus, $xb_i'\in \Fix(\varphi)\cap F'b_i'$. Now, let $x\in\Fix(\varphi)\cap F'b_i'$. Then, $x(b_i')^{-1}\in F'$ and $$(x(b_i')^{-1})\psi=x\varphi(b_i')^{-1}\varphi=x((b_i')^{-1}\varphi)=x(b_i')^{-1}b_i'((b_i')^{-1}\varphi)$$
 and so $x(b_i')^{-1}\in  X_{b_i'((b_i')^{-1}\varphi)}$.
 
 If $\Fix(\varphi)\cap F'b_i'\neq \emptyset$, then there is some $x\in F'$ such that 
   $(x\varphi)(b_i'\varphi)=(xb_i')\varphi=xb_i'$ and so $x\varphi=xb_i'((b_i')^{-1}\varphi)$. Since $F'$ is fully invariant, then $xb_i'((b_i')^{-1}\varphi)\in F'$, which yields that $b_i'((b_i')^{-1}\varphi)\in F'$.

   Clearly, $$\Fix(\varphi)=\bigcup_{i=1}^m \Fix(\varphi)\cap F'b_i',$$ and so $\Fix(\varphi)$ is computable.
\qed

\begin{remark}
\label{boundrankfix}
We remark that, since $\Fix(\psi)=\Fix(\varphi)\cap F'$, then $[\Fix(\varphi):\Fix(\psi)]\leq [G:F']$, and so  $\Fix(\varphi)$ has a generating with at most $\rank(\Fix(\psi))+[G:F']$ elements.
By \cite{[IT89]}, it follows that $\rank(\Fix(\varphi))\leq \rank(F')+[G:F']$.
\end{remark}


\section{Finite Orbits}

Given a finite orbit $\Orb_\varphi(x)$, we say that $\Orb_\varphi(x) \cap\Per(\varphi)$ is the \emph{periodic part of the orbit} and $\Orb_\varphi(x)\setminus\Per(\varphi)$ is the \emph{straight part of the orbit}.

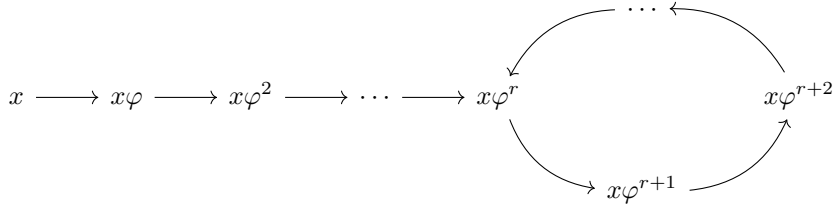
\begin{figure}[h]
\centering
  \begin{tikzcd}
    &&&&& \cdots \ar[ld,bend right] \\
    x\ar[r]
    & x\varphi \ar[r]
    & x\varphi^2 \ar[r]
    & \cdots  \ar[r]
    & x\varphi^r \ar[rd, bend right] 
    && x\varphi^{r+2} \ar[lu, bend right] \\
    &&&&& x\varphi^{r+1} \ar[ru, bend right]
  \end{tikzcd}
  \caption{A finite orbit}
\end{figure}
 
 In Figure 1, the straight part of the orbit corresponds to $\{x,x\varphi,\ldots,x\varphi^{r-1}\}$ and the periodic part of the orbit corresponds to $\{x\varphi^k\mid k\geq r\}=\{x\varphi^r,\ldots, x\varphi^{r+p-1}\}$, where $p$ is the period of $x\varphi^r$.

In \cite{[MS02]}, the authors show that, for an automorphism of a free group $F_n$, there is an orbit of cardinality $k$ if and only if there is an element of order $k$ in $\Aut(F_n)$.  Moreover, the authors prove that this result does not hold for general endomorphisms, by providing an example of an endomorphism of $F_3$ for which there is a point whose orbit has $5$ elements.
However, using a standard argument, we present a similar result for periodic parts of orbits of general endomorphisms of $F_n$. 
\begin{lemma}
Let $\varphi\in \End(F_n)$. There is a periodic point of period $k$ for some $\varphi\in \End(F_n)$ if and only if  there is an element of order $k$ in $\Aut(F_n)$. 
\end{lemma} 
\noindent\textit{Proof.}  Let $x\in \Per(\varphi)$ be a periodic point of period $k$. Consider the stable image of $\varphi$, $$S=\bigcap_{s\geq 1} F_n\varphi^s.$$
It is well known that $S$ is a free group of rank at most $n$ and that $\varphi|_S$ is an automorphism (see \cite{[IT89]}). Also, it is obvious that $\Per(\varphi)\subseteq S$, and so $x$ is a point of $S$ with a finite orbit of cardinality $k$. Therefore, by \cite[Theorem 1.1]{[MS02]}, there is an element of order $k$ in $\Aut(S)$. Since $\rank(S)=r\leq n$, then there is an automorphism of $F_n$ of order $k$, which can be defined by applying the automorphism induced by $\varphi|_S$ to the first $r$ letters and the identity in the remaining letters.

Conversely, if there is an element of order $k$ in $\Aut(F_n)$, then there is an orbit of cardinality $k$ for some automorphism of $F_n$.  Since finite orbits of automorphisms  are periodic, the result follows.
\qed\\

\begin{corollary}
\label{boundper}
There is a computable constant $k$ that bounds the size of the periodic parts of every orbit $\Orb_\varphi(x)$, when $\varphi$ runs through $\End(F_n)$ and $x$ runs through $F_n$.
\end{corollary}
\noindent\textit{Proof.}
 By \cite{[McC80]} and \cite{[Khr85]}, $\Aut(F_n)$ has an element of order $m=p_1^{\alpha_1}\cdots p_s^{\alpha_s}\in \N$, where $p_i'$s are 
 different primes, if and only if $\sum_{i=1}^s(p_i^{\alpha_i}-p_i^{\alpha_i-1})\leq n.$  We have that $\sum_{i=1}^s(p_i^{\alpha_i}-p_i^{\alpha_i-1})=\sum_{i=1}^s(p_i-1)p_i^{\alpha_i-1}$ and so, if a natural number $m\in \N$ is the order of some automorphism of $F_n$, then it only admits in its factorization primes $p$ such that $p-1\leq n$ and each of them can have exponent at most $\log_p(n)+1.$ There are finitely many integers in those conditions, and so, $m$ must be bounded above by some constant $k$ that depends only on $n$. 
\qed\\

We now prove that, given an endomorphism of a finitely generated virtually free group, we can bound the size of periodic parts of finite orbits by a computable constant. This constant depends on the endomorphism unlike the one obtained in Corollary \ref{boundper} for endomorphisms of free groups.

\begin{proposition}
\label{boundpervfree}
Let $G$ be a finitely generated virtually free group and $\varphi\in \End(G)$. Then, there is a computable constant $k$ such that the infinite ascending chain $$\Fix(\varphi)\subseteq \Fix(\varphi^{2!})\subseteq \Fix(\varphi^{3!})\subseteq \cdots$$
stabilizes after $k$ steps. Equivalently, if $x\in \Evper(\varphi)$, then the periodic part of the orbit of $x$ has cardinality at most $k$.
\end{proposition}
\noindent\textit{Proof.}   Proceeding as in the proof of Theorem \ref{computvfree}, we compute a decomposition $$F=Fb_1\cup\cdots \cup Fb_m,$$ where $F$ is a fully invariant free subgroup of $G$. We want to compute $k$ such that the ascending chain $\mc C$ defined by
$$\Fix(\varphi)\subseteq \Fix(\varphi^{2!})\subseteq  \Fix(\varphi^{3!})\subseteq \cdots$$ 
stabilizes after at most $k$ steps. For $i\in [m]$, consider the  chains $\mc C_i$ given by
$$\Fix(\varphi)\cap Fb_i\subseteq \Fix(\varphi^{2!})\cap Fb_i \subseteq \Fix(\varphi^{3!})\cap Fb_i \subseteq \cdots$$ 
Since, for all $j\in \N$, we have that $$\Fix(\varphi^{j!})=\bigcup_{i\in [m]}   (\Fix(\varphi^{j!})\cap Fb_i),$$ it follows that $\mc C$ stabilizes after $n$ steps if and only if all chains $\mc C_i$ stabilize after at most $n$ steps. 

We will prove that, for all $i\in [m]$, we can compute a constant $k_i$ such that the chain $\mc C_i$ stabilizes after $k_i$ steps and so, taking $k=\max \{k_i\mid i\in [m]\}$ suffices. 

Let $i\in [m]$. 
Since $F$ is fully invariant, we have that, for all $k\in \N$,  $(Fb_i)\varphi^k\subseteq F(b_i\varphi^k)$. Hence, the mapping $\theta: G/F\to G/F$ defined by $Fb_i\mapsto F(b_i\varphi)$ is a well-defined endomorphism. Since $G/F$ is finite, we can compute the orbit $\Orb_\theta(Fb_i)$ of $Fb_i$ through $\theta$. In particular, we can check if $Fb_i$ is periodic. If it is not, then, $\Fix(\varphi^k)\cap Fb_i=\emptyset$, for all $k\in \N$. Indeed, if there were some $k\in \N$, $x\in F$ such that $(xb_i)\varphi^k=xb_i$, then, since $x\varphi^k\in F$, we have that 
$$(Fb_i)\theta^k=F(b_i\varphi^k)=F((xb_i)\varphi^k)=F(xb_i)=Fb_i.$$
If $Fb_i$ is periodic, then let $p$ be its period and take $z\in F$ such that $b_i\varphi^p=zb_i$. Clearly, if $j\in \N$ is such that $\Fix(\varphi^j)\cap Fb_i\neq \emptyset$, then $p$ divides $j$. Also, let $C$ be the bound given  by Corollary \ref{boundper} for $n=\rank(F)+1$. Let $c$ be a letter not belonging to the alphabet of $F$ and  $\psi:F*\langle c|\rangle\to F*\langle c|\rangle$  be defined by mapping the letters of the alphabet of $F$ through $\varphi^p$ and $c$ to $zc$. 
Notice that, for all $j\in \N$, 
 $$c\psi^j=\left(\prod_{s=0}^{j-1}z\varphi^{(j-1-s)p}\right)c \quad \text{ and } \quad b_i\varphi^{jp}=\prod_{s=0}^{j-1}z\varphi^{(j-1-s)p}b_i.$$
 
 We claim that, for $x\in F$ and $q\in \N$, $$(xc)\psi^q=xc \Leftrightarrow (xb_i)\varphi^{qp}=xb_i.$$ 
 Indeed, let $x\in F$ and $q\in \N$ be such that $$ xc=(xc)\psi^q=x\varphi^{qp} c\psi^q=x\varphi^{qp}\left(\prod_{s=0}^{q-1}z\varphi^{(q-1-s)p}\right)c.$$ 
Then, $$(xb_i)\varphi^{qp}=x\varphi^{qp}\left(\prod_{s=0}^{q-1}z\varphi^{(q-1-s)p}\right)b_i=xb_i.$$

The converse is analogous.

Since, by Corollary \ref{boundper}, the periods by the action of $\psi$ are bounded above by $C$, then the periods of points in $Fb_i$ by the action of $\varphi$ are bounded above by $Cp$, which is computable since both $C$ and $p$ are. Hence, the chain $\mc C_i$ stabilizes after at most $Cp$ steps.
\qed

\begin{corollary}
\label{evperfix}
Let $G$ be a finitely generated virtually free group and $\varphi\in \End(G)$. There is a computable constant $k\in \N$ such that $\Evper(\varphi)=\Evfix(\varphi^{k!})$. 
\end{corollary}
\noindent\textit{Proof.}
Let $k$ be the constant given by Proposition \ref{boundpervfree}.  It is obvious that  $\Evfix(\varphi^{k!})\subseteq \Evper(\varphi)$. Now, let $x\in \Evper(\varphi)$. We have that there is some $s\in \N$ such that $x\varphi^s\in \Per(\varphi)$. By Proposition \ref{boundpervfree},   the period of $x\varphi^s$ is bounded above by $k$, and so it divides $k!$. Take $n\in \N$ such that $nk!>s$. This way, we have that $x\varphi^{nk!}$ belongs to the periodic part of the orbit of $x$. Thus, $x\varphi^{nk!}\varphi^{k!}=x\varphi^{nk!}$ and so, $x\in \Evfix(\varphi^{k!})$. \qed\\

Now we show that, for a fixed endomorphism, we can also bound the size of the straight part of the orbits by a computable constant, which, in combination with Proposition \ref{boundpervfree}, gives us a way of computing an upper bound on the cardinality of finite orbits.

\begin{proposition}
\label{boundstraightvfree}
Let $G$ be a finitely generated  virtually free group and $\varphi\in \End(G)$. Then, there is a computable constant $k$ that bounds the size of the straight part of every finite orbit.
\end{proposition}

\noindent\textit{Proof.} Compute a decomposition $$G=Fb_1\cup\cdots \cup Fb_m,$$ where $F$ is a fully invariant free subgroup of $G$ and write $\psi=\varphi|_F$. We can assume that $b_1=1$.
 For all $j\in \N$, consider the surjective mappings  $\varphi_j:\Img(\varphi^j)\to \Img(\varphi^{j+1})$ and  $\psi_j:\Img(\psi^j)\to \Img(\psi^{j+1})$ given by restricting $\varphi$. 
It suffices to prove that for some computable $k$, we have that $\varphi_k$ is injective  and this implies that the straight part of a finite orbit must contain at most $k$ elements.  
Indeed, suppose that there is some $x\in \Evper(\varphi)$ such that
 the straight part of $\Orb_\varphi(x)$ has
$r>k$ elements. Put $y=x\varphi^{r}$ and let $\pi$ be the period of $y$. Clearly, $y\in \Img(\varphi^{r})\subseteq \Img(\varphi^{k})$. Then $y=x\varphi^{r-1}\varphi_k$ and $y=y\varphi^{\pi-1}\varphi_k.$ But
 $y\varphi^{\pi-1}\neq x\varphi^{r-1}$ since $x\varphi^{r-1}$ belongs to the straight part of the orbit and $y\varphi^{\pi-1}$ belongs to the periodic part. This contradicts the injectivity of $\varphi_k$.

So, it remains to prove that  $\Img(\varphi^k)\simeq \Img(\varphi^{k+1})$ for some computable $k$, which, by Hopfianity of $\Img(\varphi^k)$ implies that $\varphi_k$ is injective. We have that, for all $i\in \N$,  $$G\varphi^i=F\varphi^i (b_1\varphi^i)\cup\cdots \cup F\varphi^i (b_m\varphi^i),$$ 
and so $F\varphi^i$ is a finite index subgroup of $G\varphi^i$ and $[G\varphi^i:F\varphi^i]\leq [G:F]=m$. Also,
 $0\leq \rank(\Img(\psi^{i+1}))\leq \rank(\Img(\psi^{i})),$ for every $i\in \N$. 

 Now, we describe the algorithm to compute $k$. Start by computing the smallest positive integer $j_1\in \N$  such that $ \rank(\Img(\psi^{j_1+1}))= \rank(\Img(\psi^{j_1}))$. Clearly, $j_1$ is computable: we have generators for $\Img(\psi^{i})$ for every $i\in \N$ and so we can compute its rank by computing the graph rank of its Stallings automaton. If $\rank(\Img(\psi^{j_1}))=0$, then $\psi$ is a vanishing endomorphism and so $\Img(\varphi^{j_1})$ is finite. In that case, the orbits of elements in $\Img(\varphi^{j_1})$ must be finite, since $\Img(\varphi^{k})\subseteq \Img(\varphi^{j_1})$, for $k>j_1$, and so after at most $|\Img(\varphi^{j_1})|$ iterations, we must reach a periodic point. We can compute the entire orbit of all the elements in $\Img(\varphi^{j_1})$, put $M$ to be the cardinality of the largest orbit, $k=M+j_1$ and we are done. So, suppose that $\Img(\varphi^{j_1})$ is nontrivial.

 Since free groups are Hopfian, then a free group is not isomorphic
to any of its proper quotients. Thus, $\psi_{j_1}$ must be injective. If $\Img(\varphi^{j_1})\simeq \Img(\varphi^{j_1+1})$, then, we are done. If not, by Hopfianity,  $\varphi_{j_1}$ is not injective, and so there are some $f\in F$ and $i\in [m]$ such that $((fb_i)\varphi^{j_1}\varphi)=1$. Since $\psi_{j_1}$ is injective, then $i\neq 1$. So, there is some $i\in\{2,\ldots, m\}$ such that $(fb_i)\varphi^{j_1+1}=1$ and so $b_i\varphi^{j_1+1}\in F\varphi^{j_1+1}$, thus $$[G\varphi^{j_1+1}:F\varphi^{j_1+1}]\leq m-1$$
and for all  $i\geq j_1$, $[G\varphi^i:F\varphi^i]\leq [G\varphi^{j_1+1}:F\varphi^{j_1+1}]\leq m-1$.

Now, we compute   $j_2$, the second least positive integer such that $ \rank(\Img(\psi^{j_2+1}))= \rank(\Img(\psi^{j_2}))$ and proceed as above. After $m$ steps, we have either found $k$ or we have that $G\varphi^{j_m}=F\varphi^{j_m}\simeq F\varphi^{j_m+1}=G\varphi^{j_m+1}$, and we are done.
\qed\\

Combining Proposition \ref{boundpervfree} and Proposition \ref{boundstraightvfree}, we obtain the following corollary.

\begin{corollary}
\label{upbound}
Let $G$ be a finitely generated  virtually free group and $\varphi\in \End(G)$. There is a computable constant $k$ that 
$$\max\{|\Orb_\varphi(x)| \, \big\lvert \, x\in \Evper(\varphi)\}\leq k.$$
\end{corollary}
Given an endomorphism $\varphi\in \End(G)$, we denote by $C_\varphi$ the computable constant that bounds the size of all finite orbits through $\varphi$.

\begin{corollary}
Let $\varphi\in \End(F_n)$. It is decidable whether $\Evfix(\varphi)$ is a normal subgroup of $F_n$ or not.
\end{corollary}
\noindent\textit{Proof.}
To decide if condition 1. in Proposition \ref{normalevfix} holds, we check if the generators of the free group $F_n$ are eventually fixed by computing the first $C_\varphi$ elements of their orbits. By Remark \ref{c2}, condition 2. is equivalent to $\Fix(\varphi)$ being trivial which is known to be decidable (in fact, by  \cite{[Mut21]}, we can find a basis for $\Fix(\varphi)$).
\qed

\begin{corollary}
Let $G$ be a finitely generated virtually free group and $\varphi\in \End(G)$ be defined by the image of the generators. It is decidable whether $\varphi$ is a finite order element of $\End(G)$ or not. We can also decide if $\varphi$ is aperiodic or not. 
\end{corollary}
\noindent\textit{Proof.}
The endomorphism $\varphi$ has finite order if there are $p,q\in \N$ such that for every $x\in G$, $x\varphi^p=x\varphi^q$. Since the size of finite orbits is bounded above by a computable constant, then we can check if the letters have finite orbits. If 
there is some letter $a$ with infinite orbit, then $a\varphi^p\neq a\varphi^q$, for $p,q\in \N$ with $p\neq q$. If every letter is eventualy periodic, then let $p$ be the maximum length of the straight parts of the orbits of the letters, so that $a\varphi^p$ is a periodic point for every letter $a$ and let $m$ be the least common multiple between the length of the periodic parts. Then $a\varphi^p=a\varphi^{p+m}$ for every letter and so $\varphi^p=\varphi^{p+m}$. 

Aperiodicity is similar. We have that there is an $m\in \N$ such that $\varphi^m=\varphi^{m+1}$ if and only if for each letter $a$ there is a $p$ such that $a\varphi^p=a\varphi^{p+1}$ and that is decidable simply by computing the orbits of the letters.
\qed
\begin{corollary}
Let $G$ be a finitely generated virtually free group and $\varphi\in \End(G)$. The infinite ascending chain $\Ker(\varphi)\subseteq \Ker(\varphi^2)\subseteq \ldots$ stabilizes and $\bigcup_{k=1}^{\infty} \Ker(\varphi^k)=\Ker(\varphi^{C_\varphi})$.
\end{corollary}
\noindent\textit{Proof.}
For $k>C_\varphi$, if $x\varphi^k=1$, then $x\varphi^{C_\varphi}=1$, since $|\Orb_\varphi(x)|\leq C_\varphi$.
\qed

\begin{corollary}
Let $G$ be a finitely generated virtually free group and $\varphi\in \End(G)$. Then, $\Evfix(\varphi)=(\Fix(\varphi))\varphi^{-C_\varphi}=\Ker(\varphi^{C_\varphi})\vee\Fix(\varphi)$ and $\Fix(\varphi)\simeq \faktor{\Evfix(\varphi)}{\Ker(\varphi^{C_\varphi})}$.
\end{corollary}
\noindent\textit{Proof.} We have that $\Evfix(\varphi)$ is the subgroup of points that get mapped to $\Fix(\varphi)$ by $\varphi^{C_\varphi}$. Also,  $\Ker(\varphi^{C_\varphi})\subseteq \Evfix(\varphi)$ and for every element $x\in \Evfix(\varphi)$, there is some $y\in \Fix(\varphi)$, such that $x\varphi^{C_\varphi}=y=y\varphi^{C_\varphi}$. So there must be some $z\in \Ker(\varphi^{C_\varphi})$ such that $x=yz$. Thus, $\Evfix(\varphi)=\Ker(\varphi^{C_\varphi})\vee\Fix(\varphi)$. 

Letting $\psi$ denote the restriction of $\varphi^{C_\varphi}$ to $\Evfix(\varphi)$, we have that
$$\Fix(\varphi)=\text{Im}(\psi)\simeq  \faktor{\Evfix(\varphi)}{\Ker(\psi)}=\faktor{\Evfix(\varphi)}{\Ker(\varphi^{C_\varphi})}.$$
\qed\\

We now present a result which, despite being easy and in the author's opinion, of independent interest, doesn't seem to appear in the literature. But first, we present a technical lemma.

\begin{lemma}
\label{intercosets}
Let $G$ be a group and $H,K\leq G$. Then for all $x,y\in G$, $Hx\cap Ky$ is either empty or a coset of $H\cap K$.
\end{lemma}
\noindent\textit{Proof.} Suppose that $Hx\cap Ky\neq \emptyset$ and let $z\in Hx\cap Ky$. Then, $Hz=Hx$ and $Kz=Ky$. We will prove that $(H\cap K)z=Hx\cap Ky$. Let $w\in (H\cap K)z$. Then, there are $h\in H$, $k\in K$ such that $w=hz=kz$ and so $h=k\in H\cap K$. Now, let $v\in  Hx\cap Ky$. Then $Hv=Hx=Hz$ and $Kv=Ky=Kz$ and so $vz^{-1}\in H\cap K$.
\qed

\begin{proposition}
\label{decidepreimage}
Let $G$ be a finitely generated virtually free group having a free subgroup $F$ of finite index, $H\leq G$ and $\varphi\in \End(G)$ be an endomorphism. If $\Ker(\varphi)$ is finite, then $H\varphi^{-1}$ is finitely generated. If not, the following are equivalent:
\begin{enumerate}
\item $H\varphi^{-1}$ is finitely generated
\item $H\varphi^{-1}\cap F$ is a finite index subgroup of $F$
\item  $H\varphi^{-1}\cap F$ is a finite index subgroup of $G$
\item $H\varphi^{-1}$ is a finite index subgroup of $G$
\item $H\cap G\varphi$ is a finite index subgroup of $G\varphi$
\item $H\cap F\varphi$ is a finite index subgroup of $F\varphi$
\end{enumerate}
\end{proposition}
\noindent\textit{Proof.}  It is easy to see that if $\Ker(\varphi)$ is finite, then $H\varphi^{-1}$ is finitely generated, since it is generated by the preimages of the generators of $H$ together with the kernel.

So, assume that $\Ker(\varphi)$ is infinite.
It is obvious that  $2.\Rightarrow 3.$ and that $3.\Rightarrow 4.$. It is also well known that $4.\Rightarrow 1.$ We will prove that $1.\Rightarrow 2.$, $4.\Leftrightarrow 5.$ and $5.\Leftrightarrow 6.$ and that suffices.

We start by proving that $1.\Rightarrow 2.$  Suppose that $H\varphi^{-1}$ is finitely generated. Since virtually free groups are Howson (free groups are Howson \cite{[How54]} and it is easy to see that the Howson property is preserved by taking finite extensions), then  $H\varphi^{-1}\cap F$ is also finitely generated. Then, by  Marshall Hall's Theorem, there is some finite index subgroup $H'$ of $F$ such that $H\varphi^{-1}\cap F$ is a free factor of $H'$. Take $H''$ such that $$H'=  (H\varphi^{-1}\cap F)*H''.$$ Since $H'$ is a finite index subgroup of $F$, it is finitely generated. Let $A=\{a_1,\ldots, a_n\}$ be a basis of $H'$ such that $ H\varphi^{-1}\cap F=\langle a_1,\ldots, a_k\rangle$ and $H''=\langle a_{k+1},\ldots, a_n \rangle$. 
 If $H''$ is trivial, then $H\varphi^{-1}\cap F=H'$ is a finite index subgroup of $F$ and we are done. Suppose then that $H''$ is nontrivial.  Obviously, since $\varphi$ is noninjective, $\{1\}\neq \Ker(\varphi)\trianglelefteq  H\varphi^{-1}$. Thus, $$(\Ker(\varphi)\cap F)\trianglelefteq   (H\varphi^{-1}\cap F).$$ Moreover, $\Ker(\varphi)\cap F$ is not trivial. Indeed, the fact that the kernel is infinite implies that $\varphi|_F$ is noninjective.
 Let $1\neq x\in H''$ and $1\neq y\in \Ker(\varphi)\cap F$. Then $$xyx^{-1}\in (\Ker(\varphi)\cap F)\trianglelefteq  (H\varphi^{-1}\cap F),$$ which is absurd, since the letters in $x$ don't belong to the basis set of $(H\varphi^{-1}\cap F)$.

Now, we prove that $4.\Rightarrow 5.$ Suppose that $H\varphi^{-1}$ is a finite index subgroup of $G$. Then there are  $b_i\in G$, $i\in \{0,\ldots, k\}$, such that $$G=b_0(H\varphi^{-1})\cup\cdots \cup b_k(H\varphi^{-1}). $$ 
Clearly, $$G\varphi=(b_0\varphi) (H\cap G\varphi)\cup\cdots\cup( b_k\varphi) (H\cap G\varphi)$$ and so $H\cap G\varphi$ is a finite index subgroup of $G\varphi$.

Similarly, if $H\cap G\varphi$ is a finite index subgroup of $G\varphi$, then  there are $b_i\in G$, $i\in \{0,\ldots, k\}$, such that $$G\varphi=(b_0\varphi) (H\cap G\varphi)\cup\cdots( b_k\varphi) (H\cap G\varphi).$$ Hence,  $$G=b_0(H\varphi^{-1})\cup\cdots b_k(H\varphi^{-1}),$$ because, given $x\in G$, we have that there are some $y\in G$ and $i\in \{0,\ldots, k\}$ such that $y\varphi\in H$ and $x\varphi=(b_i\varphi) (y\varphi)$ and so $x=b_iyk$ for some $k\in \Ker(\varphi)\trianglelefteq H\varphi^{-1}$. So, $yk\in H\varphi^{-1}$. Hence $5. \Rightarrow 4.$

Finally, we prove that $5.\Leftrightarrow 6.$ Assume that $H\cap G\varphi$ is a finite index subgroup of $G\varphi$. Then, there are $m\in \N$ and elements $b_i\in G\varphi$ such that $$G\varphi=(H\cap G\varphi)b_1\cup \cdots \cup (H\cap G\varphi)b_m.$$ Thus, $$F\varphi=G\varphi\cap F\varphi=\left((H\cap G\varphi)b_1 \cap F\varphi\right)\cup \cdots \cup \left((H\cap G\varphi)b_m\cap F\varphi\right).$$
By Lemma \ref{intercosets}, $H\cap F\varphi=H\cap G\varphi\cap F\varphi$ is a finite index subgroup of $F\varphi$ of index at most $m$.
Conversely, since $F\varphi$ is a finite index subgroup of $G\varphi$, then, if $H\cap F\varphi$ is a finite index subgroup of $F\varphi$, it must also be a finite index subgroup of $G\varphi$, by transitivity. Since $H\cap G\varphi\geq H\cap F\varphi$, then $H\cap G\varphi$ is also a finite index subgroup of $G\varphi$.
\qed\\

Notice that noninjective endomorphisms of free groups satisfy the hypothesis of Proposition \ref{decidepreimage}. Also, the following corollary follows directly from the proof of Proposition 
\ref{decidepreimage}.

\begin{corollary}
\label{sameindex}
The index of the subgroups in conditions 4. and 5. of Proposition \ref{decidepreimage} must coincide.
\end{corollary}

\begin{corollary}
Let $G$ be a finitely generated virtually free group and $\varphi\in \End(G)$. It is decidable whether $\Evfix(\varphi)$ (resp. $\Evper(\varphi)$) is finitely generated and, in case the answer is affirmative, a set of generators can be effectively computed.
\end{corollary}
\noindent\textit{Proof.} 
Compute a decomposition $$G=Fb_1\cup\cdots \cup Fb_m,$$ where $F$ is a fully invariant free subgroup of $G$ and write $\psi=\varphi|_F$. We can assume that $b_1=1$.

Since $F$ is fully invariant, we have that 
\begin{align}
\label{decgc}
G\varphi^{C_\varphi}=F\varphi^{C_\varphi}(b_1\varphi^{C_\varphi})\cup\cdots \cup F\varphi^{C_\varphi}(b_m\varphi^{C_\varphi}).
\end{align}

Now, it is easy to see that $\varphi^{C_\varphi}$ has finite kernel if and only if $\psi^{C_{\varphi}}$ is injective, which is decidable, since, by Hopfianity, $\psi^{C_{\varphi}}$ is injective if and only if $\rank(F)=\rank(F\psi^{C_{\varphi}})$. If $\varphi^{C_\varphi}$ has finite kernel, then it must be computable. Indeed, if, for $x\in F$, $(xb_i)\varphi^{C_\varphi}=1$, then $x\varphi^{C_\varphi}=x\psi^{C_\varphi}=b_i^{-1}\varphi^{C_\varphi}$, and so for all $i\in[m]$, we check if $b_i^{-1}\varphi^{C_\varphi}\in \Img(\psi^{C_\varphi})$, and if it is, we compute $x\in F$ such that $x\psi^{C_\varphi}=b_i^{-1}\varphi^{C_\varphi}$.
 We then have that $\Evfix(\varphi)=(\Fix(\varphi))\varphi^{-C_\varphi}$ is finitely generated and, since $\Fix(\varphi)$ is computable by Theorem \ref{computvfree} and $\Ker(\varphi)$ is computable, a set of generators for $\Evfix(\varphi)$ can be computed.  So, suppose that the kernel of  $\varphi^{C_\varphi}$ is infinite.

We know that $\Evfix(\varphi)=(\Fix(\varphi))\varphi^{-C_\varphi}$ and so, by Proposition \ref{decidepreimage}, it is finitely generated if and only if 
$\Fix(\varphi)\cap F\varphi^{C_\varphi}$ is a finite index subgroup of $F\varphi^{C_\varphi}$. Using Theorem \ref{computvfree}, we can compute a basis for $\Fix(\varphi)$, 
and so compute a set of generators for $\Fix(\varphi)\cap F\varphi^{C_\varphi}$. Notice that, since $F$ is fully invariant then $F\varphi^{C_\varphi}$ (and so $\Fix(\varphi)\cap F\varphi^{C_\varphi}$) is a subgroup of $F$. Then, we can decide if $\Fix(\varphi)\cap F\varphi^{C_\varphi}$ has finite index on $F\varphi^{C_\varphi}$
and compute right coset representatives $b_i'\in F\varphi^{C_\varphi}$
\begin{align}
\label{decfc}
F\varphi^{C_\varphi}=\left(\Fix(\varphi)\cap F\varphi^{C_\varphi}\right)b_1'\cup \cdots \cup \left(\Fix(\varphi)\cap F\varphi^{C_\varphi}\right)b_k'.\end{align}

Combining (\ref{decgc}) and (\ref{decfc}), we have that 
 $$G\varphi^{C_\varphi}=\bigcup_{i=1}^m\bigcup_{j=1}^k (\Fix(\varphi)\cap F\varphi^{C_\varphi})b_j'(b_i\varphi^{C_\varphi}),$$
 and so
$$G\varphi^{C_\varphi}=\bigcup_{i=1}^m\bigcup_{j=1}^k \Fix(\varphi)b_j'(b_i\varphi^{C_\varphi}).$$
By testing membership in $\Fix(\varphi)\cap G\varphi^{C_\varphi}$, we can check whether any two cosets coincide and refine the decomposition to obtain a proper subdecomposition where all cosets are distinct of the form (eventually relabeling the coset representatives)
$$G\varphi^{C_\varphi}=\bigcup_{i=1}^{m'}\bigcup_{j=1}^{k'} (\Fix(\varphi)\cap G\varphi^{C_\varphi})b_j'(b_i\varphi^{C_\varphi})=\bigcup_{i=1}^{m'}\bigcup_{j=1}^{k'} (\Fix(\varphi))b_j'(b_i\varphi^{C_\varphi}),$$
where $k'\leq k$ and $m'\leq m$.

For all $(i,j)\in [m']\times [k']$, we can compute $b_{i,j}\in G$ such that $b_{i,j}\varphi^{C_\varphi}=b_j'(b_i\varphi^{C_\varphi}),$ and so, by the proof of Proposition \ref{decidepreimage}, 
$$G=\bigcup_{i=1}^{m'}\bigcup_{j=1}^{k'} (\Fix(\varphi)\varphi^{-C_\varphi})b_{i,j}.$$
Having the coset representative elements $b_{i,j}$ and being able to check membership in $\Fix(\varphi)\varphi^{-C_\varphi}$, we can compute a set of generators for $\Fix(\varphi)\varphi^{-C_\varphi}=\Evfix(\varphi)$. By Corollary \ref{evperfix}, it is clear that the result also holds for $\Evper(\varphi)$.
\qed\\

Finally, we will prove that in the cases where  $\Evfix(\varphi)$ is finitely generated, we can bound its rank.  The rank of a finitely generated virtually free group $G$ is defined as the minimal cardinality of a set of generators of $G$.
\begin{proposition}
\label{boundrank}
Let $G$ be a finitely generated virtually free group, $\varphi\in \End(G)$ and $F$ be a fully invariant free subgroup of $G$.  If $\Evfix(\varphi)$ is finitely generated, then 
$\rank(\Evfix(\varphi))\leq [G:F] + \max\{\rank(F),  \rank(F)^2-3\rank(F)+3 \}$.
\end{proposition}
\noindent\textit{Proof.} If $\varphi$ is injective, then $\rank(\Evfix(\varphi))=\rank(\Fix(\varphi))\leq \rank(F)+[G:F]$, by Remark \ref{boundrankfix}. So, assume that $\varphi$ is not injective and put $\psi=\varphi|_F$. We have that $\Evfix(\psi)=\Evfix(\varphi)\cap F$ and $C_\psi\leq C_\varphi$. Thus, by \cite[Proposition 2.1]{[ASS15]}, $$[\Evfix(\varphi):\Evfix(\psi)]\leq [G:F].$$

Hence, $\Evfix(\psi)$ is finitely generated and $\rank(\Evfix(\varphi))\leq \rank(\Evfix(\psi))+[G:F].$

 If $\psi$ is injective, then $\rank(\Evfix(\psi))=\rank(\Fix(\psi))\leq \rank(F)$. So, assume that $\psi$ is noninjective. 

If $\rank(F\varphi^{C_\psi})=0$, then $\psi$ is a vanishing endomorphism, and so $\Evfix(\psi)=F$ and, in this case, we have that $\rank(\Evfix(\varphi))\leq \rank(F)+[G:F].$

If $\rank(F\varphi^{C_\psi})=1$, then $F\varphi^{C_\psi}$ is abelian, thus $\Fix(\psi)\trianglelefteq F\varphi^{C_\psi}$, and so $$\Evfix(\psi)=\Fix(\psi)\psi^{-C_\psi}\trianglelefteq F.$$
By Proposition \ref{normalevfix}, we must have that either $\Evfix(\psi)=F$, in which case we are done, or $\Evfix(\psi)=\bigcup_{k=1}^\infty \Ker(\psi^k)=\Ker(\psi^{C_\psi})$.  But we know that $\Evfix(\psi)$ is a finite index subgroup of $F$, because $\Evfix(\psi)=(\Fix(\psi))\psi^{-{C_\psi}}$ is finitely generated (see condition 4. of Proposition \ref{decidepreimage} with $G=F$). This means that $\Img(\psi^{C_\psi})$ is finite, which implies that $\psi^{C_\psi}$ is trivial. Hence, in this case $\psi$ is a vanishing endomorphism, and so  $\Evfix(\psi)=F$. In this case, we have that $\rank(\Evfix(\varphi))\leq \rank(F)+[G:F].$

Now, suppose that $\rank(F\varphi^{C_\psi})>1$. Since  $\Evfix(\psi)$ is finitely generated, then, by condition 5. of Proposition \ref{decidepreimage} (with $G=F$),  $\Fix(\psi)$ must be a finite index subgroup of $F\psi^{C_\psi}$. 
By  \cite[Proposition 3.9]{[LS77]} and \cite[Corollary 2]{[Tur96]} , we have that $$[F\varphi^{C_\psi}:\Fix(\psi)]=\frac{\rank(\Fix(\psi))-1}{\rank(F\varphi^{C_\psi})-1}\leq \rank(F)-2$$
By Corollary \ref{sameindex}, we have that $[F:\Evfix(\psi)]= [F\varphi^{C_\psi}:\Fix(\psi)]\leq \rank(F)-2 $.
But now, using \cite[Proposition 3.9]{[LS77]} again, we get that $$\frac{\rank(\Evfix(\psi))-1}{\rank(F)-1}=[F:\Evfix(\psi)]\leq \rank(F)-2,$$
and so  $\rank(\Evfix(\psi))\leq \rank(F)^2-3\rank(F)+3 $ and $$\rank(\Evfix(\varphi))\leq \rank(F)^2-3\rank(F)+3 +[G:F]. $$
\qed


\section{Further work}

The main question left open by this work concerns the generalization of these results for other classes of groups. 

\begin{problem}
Let $G$ be a  (torsion-free) hyperbolic group and $\varphi\in \End(G)$. Is it decidable whether $\Evfix(\varphi)$ (resp. $\Evper(\varphi)$) is finitely generated and, in case the answer is affirmative, can a set of generators be effectively computed?
\end{problem}

Another potentially interesting problem could be checking the existence of a better bound to the rank of $\Evfix(\varphi)$ than the one provided by Proposition \ref{boundrank}.
\begin{problem}
Is the bound given by Proposition \ref{boundrank} sharp? 
\end{problem}

\section*{Acknowledgements}
The author is grateful to Pedro Silva for fruitful discussions of these topics, which  improved the paper.

The author was  supported by the grant SFRH/BD/145313/2019 funded by Funda\c c\~ao para a Ci\^encia e a Tecnologia (FCT).%

\end{document}